\documentclass{gtart}
\usepackage{amssymb,euscript}

\input gtoutput
\volumenumber{3}\papernumber{9}\volumeyear{1999}
\pagenumbers{211}{233}\published{14 July 1999}
\proposed{Ronald Stern}
\seconded{Frances Kirwan, Walter Neumann}
\received{4 May 1999}\revised{10 June 1999}
\accepted{8 July 1999}

\newcommand{\pp}{\mathbb{P}}

\newcommand{\cc}{\mathbb{C}}
\newcommand{\rr}{\mathbb{R}}

\newcommand{\zz}{\mathbb{Z}}
\newcommand{\sss}{\mathbb{S}}

\newcommand{\ff}{\mathbb{F}}

\DeclareFontFamily{OT1}{rsfs}{}
\DeclareFontShape{OT1}{rsfs}{n}{it}{<-> rsfs10}{}
\DeclareMathAlphabet{\curly}{OT1}{rsfs}{n}{it}

\newcommand{\cdbar}{\overline{\partial}}

\newcommand{\coker}{\mathrm{coker \,}}

\newcommand{\Crit}{\mathrm{Crit}}

\newcommand{\lef}{Lefschetz }

\newcommand{\id}{\mathrm{id \,}}
\newcommand{\mod}{\mathrm{mod \,}}
\newcommand{\Fibre}{\mathrm{Fibre}}
\newcommand{\Section}{\mathrm{Section}}

\newcommand{\mgbar}{\overline{M}_g}
\newcommand{\sign}{\mathrm{sign}}

\newcounter{Universal}[section]
\renewcommand{\theUniversal}{\thesection.\arabic{Universal}}

\newenvironment{Plain}{\refstepcounter{Universal} \par \vskip 4pt
plus6.3pt minus6.3pt{\bf \theUniversal}\kern.5em}{\par \vskip7.4pt plus3pt
minus3pt}
\newenvironment{Italic}{\refstepcounter{Universal} \par \vskip 4pt
plus6.3pt minus6.3pt{\bf \theUniversal}\kern.5em\sl}{\par \vskip7.4pt plus3pt
minus3pt}  
  
\newenvironment{Thm}{\begin{Italic}{\bf Theorem}\qua\ignorespaces}
{\end{Italic}}
\newenvironment{Prop}{\begin{Italic}{\bf Proposition}\qua\ignorespaces}
{\end{Italic}}
\newenvironment{Defn}{\begin{Italic}{\bf Definition}\qua\ignorespaces}
{\end{Italic}}  
\newenvironment{Cor}{\begin{Italic}{\bf Corollary}\qua\ignorespaces}
{\end{Italic}}
\newenvironment{Lem}{\begin{Italic}{\bf Lemma}\qua\ignorespaces}
{\end{Italic}}

\newenvironment{Example}{\begin{Plain}{\bf Example}\qua\ignorespaces}
{\end{Plain}} 
\newenvironment{Rmk}{\begin{Plain}{\bf Remark}\qua\ignorespaces}
{\end{Plain}}

\newenvironment{Eqn}{\refstepcounter{Universal} $$} {\eqno \mathrm{
(\theUniversal)} $$} 

\begin{document}

\title{Lefschetz fibrations and the Hodge bundle}
\author{Ivan Smith}
\address{New College, Oxford OX1 3BN, England}
\email{smithi@maths.ox.ac.uk}

\begin{abstract}
Integral symplectic 4--manifolds may be described in terms
of Lefschetz fibrations.  In this note we give a formula
for the signature of any Lefschetz fibration in terms of the
second cohomology of the moduli space of stable curves.  As a
consequence we see that the sphere in moduli space defined by
any (not necessarily holomorphic) Lefschetz fibration has positive
``symplectic volume''; it evaluates positively with the K\"ahler
class.  Some other
applications of the signature formula and some more general results for genus
two fibrations are discussed.  
\end{abstract}

\asciiabstract{Integral symplectic 4-manifolds may be described 
in terms of Lefschetz fibrations.  In this note we give a formula for
the signature of any Lefschetz fibration in terms of the second
cohomology of the moduli space of stable curves.  As a consequence we
see that the sphere in moduli space defined by any (not necessarily
holomorphic) Lefschetz fibration has positive "symplectic volume";
it evaluates positively with the Kahler class.  Some other
applications of the signature formula and some more general results
for genus two fibrations are discussed.}



\maketitlepage


\section{Statement of Results}

The second section contains an introduction to Lefschetz fibrations
and motivation for the material of this paper, but we collect the main
results here.  Recall that all integral symplectic 4--manifolds
admit Lefschetz fibrations which in turn are equivalent to isotopy
classes of maps from a 2--sphere to the moduli space of stable curves
$\mgbar$ (satisfying appropriate conditions).  Once the genus is
sufficiently large the isotopy class of the sphere becomes a
symplectic invariant of the 4--manifold.  If the 4--manifold is
K\"ahler and the Lefschetz fibration is holomorphic the rational curve
in $\mgbar$ is a K\"ahler subvariety and in particular the
Weil-Petersson form $\omega_{WP}$ is positive on the sphere.

\begin{Thm}
For any symplectic Lefschetz fibration $f\co X \rightarrow \sss^2$
inducing $\phi_f\co \sss^2 \rightarrow \mgbar$ we have $\langle
[\omega_{WP}], [\phi_f(\sss^2)] \rangle > 0$.
\end{Thm}

\noindent The statement has geometric content; the cohomological
conditions on the sphere $\phi_f (\sss^2)$ do not alone imply the
result.  The result is a consequence of the following more general
result which extracts topological information on the 4--manifold
from the geometry of the sphere in moduli space.

\begin{Thm} \label{maintheorem}
With the notation as above, the signature of $X$ is given by
$$\sigma(X) \ = \ \langle 4c_1 (\lambda), [\phi_f(\sss^2)] \rangle -
\delta$$
where $\lambda \rightarrow \mgbar$ denotes the Hodge bundle
and $\delta$ is the number of critical fibres of the fibration.
\end{Thm}

\noindent A sketch of the proof can be found at the end of section
two.  This formula is a generalisation of one due to Atiyah for
smooth fibrations, and related work by Meyer, Wolpert and others is
well known\footnote{The formula remains valid over an arbitrary base curve
$B$ though in the sequel we shall usually leave the modifications to
the reader.}. 
In the algebraic case the theorem follows from standard Chern
class equalities;  nonetheless I have not seen the particular formula either
applied or appear in the literature; and the extension to general
symplectic fibrations (whilst unsurprising) seems new. 

The signature theorem has various implications for the geometry of
Lefschetz fibrations:

\begin{Cor}
\begin{itemize}
\item
There are no Lefschetz fibrations with monodromy group contained in
the Torelli group.

\item
Let $X \rightarrow \pp^1$ be a genus two fibration with $n=10m$
non-separating vanishing cycles and no others.
Then $X$ is naturally a double cover of
$\sss_{\mathrm{sgn}(m)}$, where
$\sss_{\mathrm{sgn}(m)}$ denotes $\sss^2
\times \sss^2$ if $m$ is even  and the non-trivial sphere bundle over
the sphere if $m$ is odd.
\end{itemize}
\end{Cor}

\noindent After earlier drafts of this work were distributed
alternative proofs 
of the first
corollary due to Li and Stipcisz also appeared.  The author's proof
also forms an Appendix to a preprint of Amoros, Bogomolov,
Katzarkov and Pantev \cite{ABKP} who formulated the statement as a Conjecture.
The existence of the double covers in the second corollary  has been
obtained, by different methods, by numerous authors but the
identification of the base of the cover in
terms of $m$ seems to have gone unnoticed.  (The more general
statement including reducible fibres is given in the text.  We also
give the classification of complex genus two Lefschetz fibrations
without reducible fibres (\ref{classify});  this result, due to Chakiris, was
rediscovered independently by the author, and to the best of his
knowledge there is no published proof.)

\ppar\textbf{Acknowledgements}\qua The material presented here is
taken from \cite{ivanthesis};  I am grateful to my supervisor
Simon Donaldson for conversations on these and related topics.  Versions
of most of this work have circulated informally and I apologise for
the delays and duplications in its independent appearance.


\section{Recalling Lefschetz fibrations}

We begin with a lightning review of the concepts underlying the rest
of the paper;  a more leisurely tour of this material is taken in
\cite{ivanthesis} amongst other sources.
Recall that a Lefschetz fibration of a smooth 4--manifold $X$
comprises a surjection to $\pp^1$, a submersion on the complement of
finitely many points $p_i$ in distinct fibres, at which there are local complex
coordinates (compatible with fixed global orientations on $X$ and
$\pp^1$) with respect to which the map takes the form $(z_1, z_2)
\mapsto z_1 ^2 + z_2 ^2$.  We always assume that the fibres contain
no $(-1)$--spheres (``relative minimality'') so in particular the fibre
genus is always strictly positive.  The existence of a
Lefschetz fibration 
guarantees that $X$ is a symplectic 4--manifold, whose topology is
determined by a monodromy homomorphism 
$$\rho_X\co \pi_1 (\pp^1 \backslash \{f(p_i) \}) \ \rightarrow \
\Gamma_g$$
where $\Gamma_g$ denotes the mapping class group of a closed
oriented genus $g$ surface.  The map $\rho_X$ maps the generators of
the fundamental group which encircle a single critical point once in
an anticlockwise fashion to positive Dehn twists in the mapping class
group. These Dehn twists are about \emph{vanishing cycles}; real
circles in a fixed fibre which shrink along some given paths to the
nodal points of the singular fibres. Thus the topology is completely
encoded in an algebraic piece 
of data which is a word in such twists in the mapping class group,
called a \emph{positive relation}.  We shall often refer to the values $\{
f(p_i)\}$---which are the critical values of $f$---by the set $\{
\Crit \}$.  The intrinsic symplectic form takes the shape
$\omega = \tau+ N f^* \omega_{\sss}$ where $\tau$ is a closed form
which is symplectic on the smooth fibres, and $\omega_{\sss}$ is
symplectic on the base $\sss^2 \cong \pp^1$.  The form is symplectic
for sufficiently large $N$, and this ``inflation'' of the horizontal
directions ensures that any local section (and hence multisection) of
the fibration with suitable orientation can be made symplectic.  The importance
of the concept for us comes from

\begin{Thm}{\rm(Donaldson)}\qua
Let $(X, \omega)$ be a symplectic 4--manifold for which $\omega$ is
the lift of an integral class $[h]$.  For sufficiently large integers $k$
the blow-up $X'$ of $X$ at $k^2 [h]^2$ distinct points admits a
Lefschetz fibration over $\pp^1$; each connected fibre, pushed back
down to $X$, is Poincar\'e dual to $k[h]$.
\end{Thm}

\noindent  The resulting Lefschetz fibration will be
relatively minimal.   We can always choose a compatible almost complex
structure 
on $X'$ such that
\begin{itemize}
\item the projection map to $\pp^1$ is pseudoholomorphic;
\item the structure is integrable in a sufficiently small tubular
  neighbourhood of each singular fibre.
\end{itemize}
Note that the exceptional sections are symplectic submanifolds of the
blow-up $X'$.  The fibres of the fibration ``downstairs'' in $X$
before blowing-up 
form a \emph{Lefschetz pencil} with finitely many \emph{base points}.
Once we have constructed a Lefschetz fibration on a symplectic
manifold there is a natural symplectic form already given to us,
without the existence result mentioned in the topological context
above; the given form and the constructible form are deformation equivalent.

\begin{Rmk}
Note that the canonical class of a symplectic manifold---which is
uniquely defined---is independent of scalings of the symplectic form.
\end{Rmk}

\noindent The choice of compatible almost complex
structures or metrics with a fixed symplectic form on $X'$ is
contractible.  Given
one such choice, the smooth fibres of the fibration become 
metric, hence conformal and complex manifolds, that is Riemann
surfaces of genus $g$.  We therefore induce a map $\pp^1 \backslash \{
f(p_i) \} \rightarrow M_g$ of a punctured sphere into the moduli space
of curves.  By the hypotheses of good local complex models, the
singular fibres of the Lefschetz fibration are naturally stable
curves (with a unique node in each) and the map completes to a map of
the closed sphere into the Deligne--Mumford stable compactified moduli
space $\mgbar$.  This map is then defined up to isotopy
independent of the 
choice of metric or almost complex structure.  The singular fibres
correspond to the intersections of the sphere with the
compactification divisor, and fall into two classes: irreducible
fibres, where we collapse a non-separating cycle in the Riemann
surface, and reducible fibres given by the one-point union of smooth
Riemann surfaces of smaller genera.  We shall see that the two kinds
of singular fibre often play a somewhat different role in the sequel;
each kind is counted by the intersection number of the sphere $\pp^1
\subset \overline{M}_g$ with the relevant components of the
compactification divisor.  Note also that the following four
stipulations are geometrically equivalent:

\begin{itemize}
\item the local complex coordinates at the $p_i, f(p_i)$ all match with fixed
  global orientations;
\item the monodromy homomorphism $\rho_X$ takes each of the obvious
  generators of the free group $\pi_1 (\pp^1 \backslash \{ \Crit \})$
  to a standard \emph{positive} Dehn twist;
\item the intersections of the sphere $\pp^1 \subset \overline{M}_g$
  with the compactification divisor of stable curves are all locally
  positive;
\item there is a symplectic structure on the total space $X'$ which
  restricts on each smooth fibre to a symplectic form.
\end{itemize}

\noindent Note that from the point of view of the moduli space of
  curves, holomorphic spheres in $\mgbar$ correspond to K\"ahler
  Lefschetz fibrations whilst smooth spheres give rise to more general
  symplectic 4--manifolds.  One can make sense of Lefschetz
  fibrations over an arbitrary base curve $B$ and \emph{mutatis
  mutandis} all the above comments apply.

A natural question to ask is how the algebraic
topology of a 4--manifold is encoded in a Lefschetz description.  By
van Kampen's theorem it is easy to see that once a set of paths and
vanishing cycles are chosen, the fundamental group of $X$ is just the
quotient of $\pi_1 (\Fibre)$ by the classes generated by the vanishing
cycles.  This in turn gives the same description for the first
homology group.  In fact all the homology groups of $X$ are given by a
pretty (short) chain complex essentially due to Lefschetz and presented
in modern notation in \cite{Zariski} (see in particular Mumford's
appendix to Chapter VI).  To recall this, fix a set of paths in the
base of the fibration and associated vanishing cycles $\delta_i$.  We have
$$0 \rightarrow H_1 (F) \stackrel{\phi}{\longrightarrow} \zz^r
\stackrel{\psi}{\longrightarrow} H_1(F)
\rightarrow 0$$
for $F$ a fixed smooth fibre and $r = \# \{ \Crit \}$.  We
define the Picard--Lefschetz twist map $T_j$ by 
$$T_j (a) = a + \langle a, \delta_j \rangle \delta_j;$$
this is the effect of the Dehn twist about the cycle
$\delta_j$ on the homology of the fibre $H_1 (F)$, and $\langle \cdot
\rangle$ denotes the intersection product.  Because the composite of
the monodromies in a Lefschetz fibration around a loop encircling all
the critical values is trivial, we have the relation 
\begin{Eqn} \label{definer}
T_r T_{r-1} \cdots T_1 \ = \ \id.
\end{Eqn}
The maps in the sequence are defined by
$$\phi\co  u \ \mapsto \ (\langle T_{i-1}T_{i-2}\cdots T_1 u, \delta_i
\rangle)_i$$ 
$$\psi\co  (a_i) \ \mapsto \ \sum_i a_i \delta_i.$$
The relation (\ref{definer}) ensures that the composite
$\psi \circ \phi$ is zero and hence the sequence does indeed give a
complex.  The cokernel of $\psi$ is precisely the first homology group
from the remarks above.  Moreover the middle homology of the complex
gives the group $G = H_2 (X) / \langle [\Fibre], [\Section] \rangle$.  To
see this, note that any element of $H_2 (X)$ in the complement of the
subspace spanned by fibres and sections projects to some graph in the
base $\sss^2$ whose endpoints are all critical values of the
fibration; such an element is closed if and only if it arises from a
union of vanishing discs bounding some homologically trivial cell in a
fibre.  Thus $G$ is indeed a quotient of $\ker \psi$.  Moreover, every
3--cell on $X$ defines on intersection with the fibres a 1--cell,
and hence the third homology of $X$ can be computed by sweeping
1--cycles in fibres around the manifold via the monodromy maps.  It
follows---again recalling the relations given by (\ref{definer})---both that the relations in $G$ are given by the image of $\phi$ whilst
the group $\ker \phi$ computes $H_3 (X)$.  Thus the word in vanishing
cycles leads to an easy computation of the homology---both Betti
numbers and torsion---for the manifold $X$.

Following this success, we ask for an expression for the
signature of $X$.  This is less straightforward.  In principle one can
compute the intersection matrix from the vanishing cycles via the
sequence above, but the formulae are highly unmanageable.  Ozbagci
\cite{Ozbagci} has shown that using Wall's non-additivity formula---a souped
up version of the Novikov additivity which gives the signature of a
4--manifold in terms of signatures of pieces resulting from cutting along a
three-manifold---one can find an algorithm for computing signature
from a word in vanishing cycles which can be fed to a computer.  In
this note we present a different formula which has the advantage of
being elegant and in closed form but which has the disadvantage of
starting not from a word in vanishing cycles but from a sphere in the
moduli space $\mgbar$.  Nonetheless we shall see that the formula
readily lends itself to certain applications to the topology and
geometry of Lefschetz fibrations.  Before giving the proof it will be
helpful to assemble some facts on signature cocycles and on the moduli
space of curves.  The natural order in which to recall these does not really
suit the proof of (\ref{maintheorem}) and for the reader's convenience
we give here the skeleton of the argument.

\ppar\textbf{A sketch of the proof}

\begin{enumerate}
\item As with the Mayer--Vietoris sequence in homology, signature can
  be computed from the pieces of a decomposition of a manifold; we
  will cut a Lefschetz fibration into its smooth part and
  neighbourhoods of critical fibres.

\item By the index theorem for manifolds with boundary, the signature
  of the smooth part can be
  expressed in terms of $\eta$--invariants of the boundary fibrations
  and the first Chern class of a determinant line bundle down the
  fibres.  (This is precisely the Hodge bundle.)

\item The determinant line bundle and $\eta$--invariant terms can be
  identified with a relative first Chern class of a topological line
  bundle defined by a signature cocycle in the group cohomology of the
  symplectic group.  Thus the signature of the smooth part of the
  fibration is computed by a symmetric function on the symplectic
  group whose arguments are the monodromies around the boundary circles.

\item By naturality properties of the cocycle and the above, the
  difference between evaluating the first Chern class of the Hodge
  bundle on a surface in $\mgbar$ and the relative first Chern class
  on a surface 
  with boundary in $M_g$ is entirely determined by the conjugacy
  classes of the monodromies in $Sp_{2g}(\zz)$; since these are all
  equal for a Lefschetz fibration, the discrepancy is measured by a
  single integer for each genus.

\item This integer is fixed by determining the theorem for at least
  one fibration at every genus; a Riemann--Roch theorem gives the
  theorem for all projective fibrations, completing the proof.
\end{enumerate}


\section{The signature cocycle}

Much of the material here is drawn from Atiyah's pretty discussion in
\cite{Atiyaheta}. We start with Novikov's additivity formula, which
states that if we decompose a 4--manifold $X = X_1 \cup_Y X_2$ along
an embedded 3--manifold $Y$, then
$$\sigma(X) = \sigma(X_1,Y) + \sigma(X_2, Y)$$
where the left hand side denotes the signature of $X$ and
for a manifold with boundary $(Z, \partial Z)$ the relative signature
$\sigma(Z, \partial Z)$ denotes the difference between the number of
positive and negative eigenvalues on the intersection form in the
middle cohomology of $Z$ (but this form is no longer non-degenerate).
Alternatively, on the image of the relative cohomology $H^2 (Z,
\partial Z) \hookrightarrow H^2 (Z)$ the intersection form is
non-degenerate and we take the signature of the form restricted to
this subspace.  

Let $X \rightarrow \sss^2$ be a Lefschetz fibration.  If we
decompose $X$ into tubular neighbourhoods of the various singular
fibres and a swiss cheese then by Novikov additivity we write the
signature as the signature of the fibration over $\sss^2 \backslash \{
\mathrm{Discs} \}$ corrected by the sums of the local signatures.  An
easy computation, noting that the total space
of such a tube retracts to the singular fibre, gives that

\begin{itemize}
\item $\sigma (Z, \partial Z) = -1$ for a neighbourhood of a reducible
  singular fibre (separating vanishing cycle)
\item $\sigma(Z, \partial Z) = 0$ for a neighbourhood of an
  irreducible singular fibre (non-separating vanishing cycle).
\end{itemize}

\noindent Thus we know $\sigma(X) = -s + \sigma(W, \partial W)$ where
  $s$ is the number of separating vanishing cycles and $W$ is the
  preimage $f^{-1} (\sss^2 \backslash D)$, $D$ a neighbourhood of the
  set of critical values.  Following \cite{Atiyaheta} it is now natural to
  introduce the notation $\sigma(A_1,\ldots,A_r)$ for the signature of
  the 4--manifold $Z$ which is a fibration by genus $g$ curves over a
  sphere with $r+1$ ordered open discs deleted, and for which the
  monodromies around the first $r$ boundary circles are given by
  elements $A_i \in \Gamma_g$.  It follows of course that the
  monodromy about the last boundary circle is just the inverse of the
  ordered product of these matrices, for the fundamental group of the
  base is free.  Write $\sss^2 \backslash D = \Sigma$.
  The topology of $Z$ is
determined by a representation $\pi_1 (\Sigma) \rightarrow \Gamma_g$
which may be composed with the standard representation $\Gamma_g
\rightarrow Sp_{2g}(\zz)$; this gives a flat vector bundle $E$  over
$\Sigma$ with fibre the first complex cohomology of the fibre.  On this
vector bundle there is a non-degenerate skew-symmetric
form given by the cup-product $\langle \cdot \rangle$ on $H^1 (\Fibre)$.
Combining this skew form with the cup-product on classes from the base we
  obtain an indefinite Hermitian
  structure on the single vector space $H^1(\Sigma, \partial
\Sigma; H^1(\Fibre))$ given by the cohomology with local
  coefficients in $E$.

\begin{Lem} {\rm\cite{Atiyaheta}}\qua
The signature $\sigma(A_1 , \ldots, A_r)$ of $Z$ is the signature of
the Hermitian form defined above on the cohomology group
with local coefficients $H^1 (\Sigma, \partial
\Sigma; H^1(\Fibre))$.
\end{Lem}

\noindent The proof amounts to a careful application of the Leray-Serre
spectral sequence for a fibration.  That only the homologically
non-trivial monodromies enter is unsurprising; signature after all is
a homological invariant.
It follows that the signature of a Lefschetz fibration is completely
determined by the number of separating vanishing cycles $s$ and the
value of the function $\sigma(\{ A_i \})$ where the $A_i$ can now be
taken to be the symplectic matrices corresponding to the monodromies
about homologically essential vanishing cycles.  Applying the Novikov
property to this symmetric function on the symplectic group gives the
following critical relation:
$$\sigma(A_1, A_2, A_3) \ = \ \sigma(A_1, A_2) + \sigma(A_1 A_2,
A_3)$$
and hence, splitting the sphere with four holes in two
distinct ways,
$$\sigma(A_1, A_2) + \sigma(A_1 A_2, A_3) \ = \ \sigma(A_2,A_3) +
\sigma(A_2 A_3, A_1).$$
This last formula is exactly the \emph{2--cocycle} condition
in group cohomology, and it follows that we have defined an element
$[\sigma] \in H^2 (Sp_{2g}(\rr); \zz))$.  Such group cohomology
elements correspond (cf Lemma (\ref{perfect})) to central extensions
of the group by the integers
and we have an associated sequence
$$0 \rightarrow \zz \rightarrow Sp^{\sigma} \rightarrow Sp_{2g}(\rr)
\rightarrow 0;$$
moreover again by standard properties of group cohomology
\cite{Brown} there
is a section to the last map $c^{\sigma}\co  Sp_{2g} \rightarrow
Sp^{\sigma}$ with the properties that

\begin{itemize}
\item the product $c^{\sigma} (A_1) c^{\sigma}(A_2) c^{\sigma}(A_1
  A_2)^{-1}$ gives a well-defined element of the central factor $\zz$
  for any $A_i \in Sp_{2g}(\rr)$,
\item $\sigma(A_1, A_2) \ = \ c^{\sigma} (A_1) c^{\sigma}(A_2) c^{\sigma}(A_1
  A_2)^{-1}$ for any $A_i \in Sp_{2g}(\rr)$.
\end{itemize}

\noindent The central point is that this section $c^{\sigma}$ has an
  interpretation in terms of a line bundle.  Recall that the
  cohomology groups down the fibres of $Z \rightarrow \Sigma$ formed a
  flat vector bundle $E$ with a symplectic structure.  If we
  complexify $E$ then the form $i \langle \cdot \rangle$ is Hermitian
  of type $(g,g)$ and we can choose a splitting $E_{\cc} = E^+ \oplus E^-$
  into maximal positive definite and negative definite subspaces.
  Such splittings correspond to reducing the structure group of
  $E_{\cc}$ from $U(g,g)$ to $U(g) \times U(g)$; the quotient
  homogeneous space is contractible so such splittings necessarily
  exist and are unique up to homotopy.

\begin{Lem} \label{number1}{\rm\cite{Atiyaheta}}\qua
Let $L$ be the line bundle $(\det E^+) \otimes (\det E^-)^{-1}$.
The section $c^{\sigma}$ defines a homotopy class of trivialisations
of $L^2$ over any loop $\gamma \subset \Sigma$ and hence a relative
first Chern class $c_1 (L^2; c^{\sigma}) \in H^2 (\Sigma, \partial
\Sigma) \cong \zz$.  Then working over a sphere with three discs deleted
$$\sigma (A_1, A_2) = c^{\sigma} (A_1) c^{\sigma}(A_2) c^{\sigma}(A_1
  A_2)^{-1} = c_1 (L^2; c^{\sigma})$$
for any $A_i \in Sp_{2g}$. 
\end{Lem}

\noindent Since we have Novikov additivity for $\sigma(\{ A_i \})$ and
  relative Chern classes behave well with respect to connected sums,
  it follows that for any fibration $Z \rightarrow \Sigma$ we have
\begin{Eqn} \label{identity}
\sigma (Z, \partial Z) \ = \ \sigma ( \{ A_i \}) \ = \ c_1 (L^2;
c^{\sigma}).
\end{Eqn}
Therefore to understand the signature of a Lefschetz
fibration we need only understand this relative first Chern class and
its trivialisation over loops as defined by $c^{\sigma}$.  Before
turning to this we recall one well known observation which is relevant
both above and in the sequel.  The section $c^{\sigma}$ is in fact
unique, since the difference between any two choices would give a
homomorphism from the group $Sp_{2g}$ to the integers.  But no such
homomorphisms can exist; for there is a canonical homomorphism from
the mapping class group $\Gamma_g$ onto $Sp_{2g}$ and

\begin{Lem} \label{perfect}
The mapping class group $\Gamma_g$ is perfect for $g \geq 3$ and has
finite cyclic abelianisation for $g = 1,2$.
\end{Lem}

\noindent This result is well known; $(\Gamma_1)_{ab} = \zz_{12}$
whilst $(\Gamma_2)_{ab} = \zz_{10}$.  For $g \geq 3$ the usual proof
simply writes a generating Dehn twist as an explicit product of
commutators.  More in line with our thinking is the following
geometric sketch.  It was shown by Wolpert \cite{Wolpert} that for
genus $g \geq 3$ the second cohomology of the moduli space of stable curves is
generated by a K\"ahler class and the divisors given by the components
of the compactifying locus of stable curves.  By Poincar\'e duality it
follows that there are homology elements $[C_i]$ which have algebraic
intersection $1$ with the $i$-th component of the compactification
divisor and $0$ with all others.  We can represent these classes by
embedded surfaces; since the moduli space is of high dimension (and
the orbifold loci of high codimension) we can tube away excess
intersections, at the cost of increasing genus, and even find
representing surfaces $C_i$ with geometric intersection numbers with
the stable divisor the same as the algebraic intersection numbers.
But the fundamental group of a surface of genus $g$ with one boundary
circle admits the presentation
$$\langle a_1, b_1, \ldots , a_g, b_g, \partial \ | \ \prod [a_i, b_i]
= \partial \rangle.$$
Thus in the fibration of curves over the surface $C_i$
defined by the universal property of the moduli space, the monodromy
about the unique singular fibre---a standard positive Dehn twist---is
expressed as a product of commutators.  Since we can do this for any
isotopy class of Dehn twist the result follows.


\section{Hodge lines and determinant lines}

To warm up we will treat the case of a projective fibration $f\co  X
\rightarrow B$ and
introduce some of the objects that appear in the final statement of
the main theorem.  (In fact the proof of (\ref{maintheorem}) will not
crucially rely on a
separate treatment for the projective case but some aspects are
simplified by giving such an argument, and it gives some structure to
the theory of the Hodge bundle we want to quote.)  Now $df$ has
maximal rank only
away from $\{ \Crit \}$ the set of critical points $p_1, \ldots, p_r$ of
$f$, but there is an exact sequence of \emph{sheaves}
\begin{displaymath}
0 \rightarrow f^* \Omega^{1,0}_B \rightarrow \Omega^{1,0}_X
\rightarrow \Omega^1 _{X/B} \rightarrow 0
\end{displaymath}
where the last term is defined by the sequence.  Away from
$\{ \Crit \}$ there is an isomorphism $\Omega^1 _{X/B} \cong \mathcal{K}_X
\otimes f^* T^{1,0}_B$ of bundles.

\begin{Defn}
The line bundle $\omega_{X/B} = \mathcal{K}_X \otimes f^* T^{1,0}_B$
is called the \emph{dualising sheaf} of $f$.
\label{dualis}
\end{Defn}

\noindent The adjunction formula says that if
$C$ is a smooth curve in a complex surface $Z$ then $\mathcal{K}_C  =
\mathcal{K}_Z \otimes \mathcal{O}_Z (C) |_C$.  For singular $C$ we
\emph{define} the ``dualising sheaf'' for $C$ via $\omega_C =
\mathcal{K}_Z \otimes \mathcal{O}_Z (C) |_C$ where $\mathcal{O}_Z (C)$
is still the line bundle associated to the divisor $C$. Despite the
definition this is naturally associated to $C$, independent of its
embedding in any ambient surface.  Since in a
Lefschetz fibration all the fibres have trivial normal bundles, it
follows that $\omega_{X/B} |_{X_b} = \omega_{X_b}$ for each fibre $X_b
= f^{-1}(b)$.  The results of
the following theorem are drawn from \cite{BPV}, Chapter III, sections
11 and 12.  Write $f_* \mathcal{F} = R^0 f_* \mathcal{F}$.

\begin{Thm}
Let the notation be as above;  suppose $f$ has generic fibre genus $g$.
\label{results}
\begin{enumerate}
\item $f_* \omega_{X/B}$ is locally free of rank $g$.
\item $R^1 f_* \omega_{X/B}$ is the trivial line bundle
$\mathcal{O}_B$ (and the higher direct images vanish for
reasons of dimension). \qed
\end{enumerate}
\end{Thm}

\noindent  These facts require that $f$ have connected
fibres.  The following definition is formulated under the naive
assumption that the universal curve $\mathcal{C} M_g \rightarrow M_g$
exists and extends to the compactification.  Whilst this fails because
of orbifold problems, the failure is in a real sense only technical; we could
work instead with specific families and the ``moduli functor''
\cite{arabel} if necessary.

\begin{Defn}
\label{hodge}
Write $\lambda_B = \wedge^g f_* \omega_{X/B}$, the \emph{Hodge bundle}
over $B$.  If the fibration $f$ gives a map $\phi_f\co  B \rightarrow
\overline{M}_g$ then $\lambda_B = \phi_f^* \lambda$ for $\lambda$
the Hodge bundle of the universal curve
$\mathcal{C}\overline{M}_g \rightarrow \overline{M}_g$.
\end{Defn}

\noindent The Hodge bundle is well known to generate the Picard group of line
bundles on $M_g$ and to extend to $\overline{M}_g$.  We
have two cohomology classes $c_1(\lambda), \delta \in H^2
(\overline{M}_g)$ where $\delta$ ``counts intersections with the
divisor of stable curves''; that is, $\delta$ is the first Chern class
of the line bundle with divisor $\overline{M}_g \backslash M_g$, or
the Poincar\'e dual of the fundamental class of $\overline{M}_g
\backslash M_g$.
For projective fibrations Mumford deduces the following from the
Grothendieck--Riemann--Roch theorem \cite{MumfordBirk}:

\begin{Prop}
Let $f\co  X \rightarrow B$ be a complex \lef fibration of a projective
surface.  Let $[\Crit]$ denote the cohomology class defined by the
critical points of $f$ viewed as a subvariety of $X$.  Then
\begin{displaymath}
12 c_1 (\lambda_B ) \ = \ f_* \big ( c_1 ^2 (\omega_{X/B}) + [\Crit]
\big ).
\end{displaymath}
\label{stage1}
\end{Prop}\vglue -0.5cm

\noindent The geometry here enters in the form of an exact sequence
\begin{Eqn} \label{exactseq}
0 \rightarrow \Omega^1 _{X/B} \rightarrow \omega_{X/B} \rightarrow
\omega_{X/B} \otimes \mathcal{O}_{\{ \Crit \}} \rightarrow 0.
\end{Eqn}
For the structure sequence of $\{ \Crit \} \subset X$
gives 
\begin{displaymath}
0 \rightarrow \mathcal{I}_{\{ \Crit \}} \rightarrow \mathcal{O}_X \rightarrow
\mathcal{O}_{\{ \Crit \}} \rightarrow 0
\end{displaymath}
\noindent where $\mathcal{I}_{\{ \Crit \}}$ is the ideal sheaf of $\{ \Crit \}$ in
$X$.  The result follows on tensoring with $\omega_{X/B}$ provided we
show that $\mathcal{I}_{\{ \Crit \}} \otimes \omega_{X/B} = \Omega^1 _{X/B}$.  Now
away from $\{ \Crit \}$ the sequence
\begin{displaymath}
0 \rightarrow f^* \Omega^{1,0} _B \rightarrow \Omega^{1,0}_X
\rightarrow \Omega^1 _{X/B} \rightarrow 0
\end{displaymath}
\noindent is an exact sequence of \emph{bundles} and $\Omega^1 _{X/B},
\ \omega_{X/B}$ coincide.  Near a point of $\{ \Crit \}$ $f$ has local form
$(z_1 , z_2 ) \mapsto z_1 z_2 $ for suitable coordinates.
Accordingly 
\begin{displaymath}
\Omega^1 _{X/B} \ = \ \coker \ \big ( f^* \Omega^{1,0}_B \rightarrow
\Omega^{1,0}_X \big ) \ = \ \frac{\mathcal{O}_X dz_1 + \mathcal{O}_X
dz_2 }{\mathcal{O}_X . (z_1 dz_2 + z_2 dz_1 )}.
\end{displaymath}
Now $\mathcal{I}_{\{ \Crit \}}$ is by definition the sheaf of
ideals generated by $\langle z_1 , z_2 \rangle $ over $\mathcal{O}_X$,
whilst the 
dualising sheaf $\omega_{X/B} = \mathcal{O}_X . (dz_1 \wedge dz_2 )
\otimes f^* ( \partial / \partial t)$ for $t$ a coordinate on $B$.
Explicitly
\begin{displaymath}
f^* \frac{\partial}{\partial t} = \frac{1}{2z_2
}\frac{\partial}{\partial z_1 } + \frac{1}{2z_1
}\frac{\partial}{\partial z_2 } \ \Longrightarrow \ \omega_{X/B} =
\mathcal{O}_X . \left ( \frac{dz_2 }{z_2 } - \frac{dz_1 }{z_1 } \right
).
\end{displaymath}
Now consider the transformation on $\mathcal{O}_X \cdot
\langle dz_1 ,
dz_2 \rangle $ generated by
\begin{displaymath}
dz_1 \mapsto z_1 \frac{dz_2 }{z_2 } - dz_1 \ \ ; \ \ dz_2 \mapsto z_2
\frac{dz_1 }{z_1 } - dz_2
\end{displaymath}
and see the local forms for $\Omega^1 _{X/B}$ and
$\mathcal{I}_{\{ \Crit \}} \otimes  \omega_{X/B}$ are clearly equal. 

The Grothendieck--Riemann--Roch theorem states that
\begin{displaymath}
ch(f_! \mathcal{F}) \ = \ f_* \big ( ch \mathcal{F} . \mathcal{T} (T^1
_{X/B}) \big ) \qquad \qquad \textrm{[GRR]}
\end{displaymath}
for $\mathcal{F}$ any coherent sheaf on an irreducible
non-singular projective variety, $ch$ the Chern character and
$\mathcal{T}$ the Todd class; recall also that $f_! \mathcal{F} = \sum
(-1)^i R^i f_* \mathcal{F}$.  The relative tangent sheaf $T^1 _{X/B}$
of the map $f$ is defined via $T^1 _{X/B} =  T^{1,0}_X -
f^* T^{1,0}_B$ as an element of $K$--theory.  The usual expansions of
Chern characters and Todd classes, combined with taking total Chern
classes in the exact sequence (\ref{exactseq}), give Mumford's result.

\begin{Cor}
The main result (\ref{maintheorem}) is valid for projective
fibrations.
\end{Cor}

\begin{proof}  By Hirzebruch's classical theorem we know 
\begin{displaymath}
\sigma (X) \ = \ \frac{\langle p_1 (TX), [X] \rangle }{3}.
\end{displaymath}
Using $p_1 (TX) = c_1 ^2 (T^{1,0}_X) - 2c_2 (T^{1,0}_X)$ and
from the definition of 
$\omega_{X/B}$ it follows that
\begin{eqnarray*}
\ & c_1 (T^{1,0}_X ) \ = \ f^* c_1 (T^{1,0}_B ) - c_1 (\omega_{X/B})
\\ 
\Longrightarrow & \  c_1 ^2 (\omega_{X/B}) \ = \ p_1 (TX) +
2c_2 (T^{1,0}_X ) + 2f^* c_1 (T^{1,0}_B ) . c_1 (\omega_{X/B}) \\
\Longrightarrow & \ f_* c_1 ^2 (\omega_{X/B}) \ = \ f_* p_1 (TX) +
2f_* c_2 (T^{1,0}_X ) + 2c_1 (T^{1,0}_B ). f_* c_1 (\omega_{X/B}).
\end{eqnarray*}
Now by the remarks after (\ref{dualis}), $\omega_{X/B} |_{X_b}
= \omega_{X_b} = \mathcal{K}_X |_{X_b}$ and since the canonical
divisor of a genus $g$ curve has degree $2g-2$ it follows that 
\begin{displaymath}
f_* c_1 (\omega_{X/B}) \ = \ (2g-2)[B] \ \in \ H^0 (B, \mathbb{Q}).
\end{displaymath}
Moreover in the top dimension $f_*$ commutes with evaluation
and recalling $\chi(X) = \chi(B) \chi(F) + \delta$ the result follows.
\end{proof}

\noindent Projective fibrations exist with fibres of every
genus $g > 0$ so we have now proven the result for at least one
fibration of each genus.  We now turn to the general smooth case in
the framework of differential geometry.
The signature of any Riemannian manifold can be defined as the index of a
differential operator $d+d^*\co  \Omega^+ \rightarrow \Omega^-$, where
$\Omega^{\pm}$ are the eigenspaces for an involution $\phi \mapsto i^{\flat}
\phi$ (some suitable power $\flat = \flat(p)$) and $\phi \in
\Omega^p$.  On a complex Riemann surface $M$ this signature operator
is equivalent to
$$\cdbar\co   \Omega^{0,0} (M) \oplus \Omega^{1,0}(M) \ \rightarrow \
\Omega^{0,1}(M) \oplus \Omega^{1,1}(M).$$
The determinant line of this operator has fibre (using Serre
duality)
$$\mathcal{L}_{\mathrm{det}} \ = \ [\det H^0 (M, \Omega^1)]^{-2}$$
thus fibrewise there is a naturally defined isomorphism
$\lambda \ \equiv
\mathcal{L}_{\det}^{-1/2}$.

\noindent Now in general a determinant line bundle
admits a canonical metric and connexion \cite{Freed}, and the holonomy
of the connexion is given by an expression of the shape $\exp
(i\eta(\cdot))$; that is, the exponential of an $\eta$--invariant of the
boundary manifold.  Usually there is no canonical choice of logarithm
for this holonomy.  However, for the particular case of the signature
operator, such a canonical logarithm \emph{does exist}.  The reason
for this special behaviour is
that the zero-eigenforms for the relevant differential operator give
rise to harmonic forms, which by Hodge theory are governed by the
topology of the manifold; thus the dimensions of the zero-eigenspaces
cannot jump as for a general differential operator.

The upshot is a canonical trivialisation via $\eta$--invariants for the
determinant line bundle $\mathcal{L}_{\det}$ over any circle and hence
the boundary of
$\Sigma$.  The index theorem for manifolds with boundary \cite{APS}
has famously been used to give a formula for the signature in terms of
$L$--polynomials with a boundary correction term defined via
$\eta$--invariants.  Comparing the terms of this expression to the
universal expression for the first Chern form of a determinant line
bundle gives the central identity (compare to (\ref{identity})!)
\begin{Eqn}
\label{secondtime}
\sigma(Z, \partial Z) \ = \ -2c_1(\mathcal{L}_{\det}; \eta)
\end{Eqn}
where the notation refers to a relative Chern class defined
with respect to the $\eta$--invariant trivialisations over $\partial \Sigma$.

\noindent With all this preamble we can forge the bridge between the
two approaches.
The determinant line bundle
comprises a piece $\mathcal{L}'$ corresponding to non-zero eigenvalues
of the differential operator which is canonically trivialised
(topologically though not as a unitary bundle) by Quillen's canonical
determinantal section, which by construction is non-vanishing there.
Thus as a topological bundle the determinant line bundle is isomorphic
to the bundle given by taking only the zero-eigenvalue spaces (of
harmonic forms in $\Omega^{\pm}$ in our case):
\begin{Eqn} \label{number2}
\mathcal{L}_{\det} \ \cong_{\mathrm{TOP}} \ (\det H^-) \otimes (\det
H^+)^{-1}.
\end{Eqn}
\begin{Lem}
The line bundle $\mathcal{L}_{\det}$ is topologically the dual of the
line bundle $L_{\sigma}$ defined from the signature cocycle.  Moreover
the trivialisations of $\mathcal{L}_{\det}^2$ and $L_{\sigma}^{-2}$
(defined by $\eta$--invariants and $c^{\sigma}$ respectively) coincide.
\end{Lem}

\begin{proof}
That the line bundles are dual amounts to an identification of the
positive and negative harmonic forms of $\Omega^{\pm} \cap H^1
(\Fibre)_{\cc}$ with the 
positive and negative eigenspaces for the Hermitian form $i \langle
\cdot \rangle$; then compare the formulae defining $L$ in
(\ref{number1}) and (\ref{number2}).  But by definition the
$\Omega^{\pm}$ are defined to be eigenspaces for an operator whose
index is signature and the $E^{\pm}$ comprising $L$ are the definite
subspaces for the Hermitian form arising from the symplectic signature
pairing.  That the trivialisations agree is a consequence of the lemma
(\ref{perfect}).  For if the two trivialisations differed then their
difference would define a map from the set of components of the
boundary $\partial \Sigma$ to $\zz$ depending only on the particular
monodromies associated to these components.  Moreover since we know
that there are identities (\ref{identity}, \ref{secondtime}) for all
fibrations $Z \rightarrow \Sigma$ we always know that
$$\sum_{\gamma \in \pi_0 (\partial \Sigma)}
\eta_{\mathrm{triv}}(\gamma) + c^{\sigma}_{\mathrm{triv}}(\gamma) \ =
\ 0.$$
Since the values of the trivialisations on loops depend only
and naturally on the monodromies, this difference map defines $N\co 
\Gamma_g \rightarrow \zz$ which (taking the above relation for
$\Sigma$ a sphere with three discs removed) is a homomorphism.  But by
(\ref{perfect}) we
know the mapping class groups admit no such non-trivial homomorphisms.
\end{proof}

\noindent Assembling our various identifications we have proven
$$\sigma(X) \ = \ -s + 4\langle c_1(\lambda), [\Sigma] \rangle$$
where we interpret the right hand side with respect to the
$\eta$--trivialisation still (and recall $s$ is the number of
separating vanishing cycles).  
To produce the final formula (\ref{maintheorem}) we now need
to understand the exact nature of the holonomy term for one of the
Dehn twist monodromies in a Lefschetz fibration.  Atiyah gives the
precise formula which shows the sense in which the $\eta$--invariant
gives a canonical logarithm for the holonomy of the Quillen connexion:
$$\delta \log \det D_{\mathrm{Quillen}} \ = \ -i\pi \eta - \frac{1}{2}
 \log_{\sign} (\mathrm{Monod}(\curly{H}));$$
here in the final term, which is an integer, $\curly{H}$ is the
 ``topological  determinant bundle'' $\curly{H} = (\det H^+)^{-1}
 \otimes \det H^-$.  The monodromy
 denotes the particular element of $Sp_{2g}(\zz)$ corresponding to the
 fibration over a given boundary circle, and $\log_{\sign}$ denotes
 a choice of logarithm for this monodromy given by the explicit
 signature cocycle we began with.  The final answer can therefore be
 given by a direct computation with this cocycle.  More simply, given
 the work at the start of the section and the naturality properties of
 the signature of manifolds (and hence the cocycle), we know that the
 answer depends 
 only on the conjugacy class of the monodromy in the mapping class
 group.  Since all our monodromies are positive Dehn twists about
 non-separating vanishing cycles, we are interested in a single
 integer for each genus $g$.  This is then determined by the
 computation of the signature 
 for a single Lefschetz fibration with genus $g$ fibres. But we
 already know the answer for projective fibrations: writing $n$ for the number
 of non-separating vanishing cycles,
$$\sigma(X)  \ = \ -s -n + \langle 4c_1 (\lambda), [\pp^1] \rangle$$
which is just as we require.

\noindent Note that from this point of view the singular fibres enter
the formula from naturally different perspectives; the separating ones because
they affect $H^2$ and invoke a local contribution to signature, the
non-separating ones because they affect $H^1$ and hence the global
monodromy which detects the extent to which the manifold is not
homologically a product. 


\section{Applications to genus two fibrations}

In the last sections we give some applications of the signature
formula and digress into some of the topics we encounter.

\begin{Example}
The moduli space of genus two curves is in fact globally
\emph{affine}, and the ample divisor given by $f_* c_1^2 (\omega)$,
for $\omega$ the dualising sheaf of the universal curve, is empty on
$M_2$ and a linear combination of boundary divisors on the
compactification.  Mumford \cite{MumfordBirk} has calculated precisely
that at genus two
$$10 c_1(\lambda) = \delta_0 + 2 \delta_1$$
where $\delta_0, \delta_1$ are respectively the irreducible
components of $\overline{M}_2 \backslash M_2$ corresponding to curves
which are generically irreducible or a union of two elliptic curves
respectively.  From the signature formula, we therefore see that the
signature for a genus two Lefschetz fibration is determined completely
by the numbers of non-separating and separating vanishing cycles $n$
and $s$:
\begin{Eqn} \label{fracsign}
\sigma(X) \ = \ - \frac{3}{5} n - \frac{1}{5} s.
\end{Eqn}
This ``fractional signature formula'' was first established
by Matsumoto \cite{Matsumoto} by related but different means.
\end{Example}

\noindent We can now compare the two different formulae for signatures
of genus two fibrations to some effect.  Recall that $(\Gamma_2)_{ab}
= \zz_{10}$ and hence for any genus two fibration we have $10 | (n+2s)$.

\begin{Prop}
Let $(X,f)$ be a genus two Lefschetz fibration with $n+2s=10m$.  Write
$\sss_{\mathrm{sgn}(m)}$ for the product $\sss^2 \times \sss^2$ if
$m$ is even and for the non-trivial sphere bundle over the sphere if
$m$ is odd.  Then
$$X \# s \overline{\cc \pp}^2  \stackrel{2:1}{\longrightarrow}
\sss_{\mathrm{sgn}(m)} \# 2s \overline{\cc \pp}^2 \ \supset C;$$
the blow-up of $X$ at $s$ points admits a smooth double
cover over the blow-up of the relevant sphere bundle over a sphere at
$2s$ points, ramified over a smooth surface $C$.
\end{Prop}

\noindent Equivalently, $X$ is given by blowing down $(-1)$--spheres in the
fibres of a fibration arising from double covering
a non-minimal rational surface over a curve $C$ which is the canonical
resolution of singularities of a curve $C' \subset
\sss_{\mathrm{sgn}(m)}$ containing $s$ infinitely close 
triple points.  (Such a point is given by the singularity at the
origin of $z_1 ^3 + z_2 ^6 = 0$; the curve has $3$ sheets meeting
mutually tangentially at this point, and the sheets are separated by
two successive blow-ups.)
Thus if $X$ has no reducible fibres then it has
$10m$ singular fibres for some $m$, and $X$ double covers the sphere
bundle with ``parity'' the same as the parity of $m$, branched over a
smooth two-dimensional surface $C$.  As a piece of notation, recall
that $\ff_k = \pp ( \mathcal{O} \oplus \mathcal{O}(k))$ denotes the
unique complex (or symplectic) ruled surface with symplectic sections of
self-intersection $\pm k$; moreover all projective bundles over $\pp^1$
are the projectivisations of (not uniquely determined) vector bundles.

\begin{proof}
The proof of the proposition is reasonably straightforward
and versions due to Fuller \cite{Fuller} and Siebert--Tian \cite{ST} have now
appeared (an independent proof was given in \cite{ivanthesis}).  The
idea is simple; on choosing a metric  all the smooth fibres are
hyperelliptic Riemann surfaces and admit
natural double branched covers over spheres.  These can be patched
together smoothly except near separating singular fibres.  The point
is that the map fibrewise is given by sections of the canonical bundle
on the fibre, and the nodal points in reducible fibres are base points
of the canonical system on a stable curve; for smooth or stable
irreducible curves the canonical linear system has no base points and
the branch locus varies smoothly.  Near reducible fibres we can assume
the complex structure is integrable and graft in a local holomorphic
model.

\noindent This argument (and those of Fuller and Siebert--Tian) gives
the base of the double cover the structure of a sphere bundle over the
sphere but does not identify it beyond that.  The signature formulae
allow one to do precisely this. The ruled surface $\ff_k = \pp(E) = \pp
  (\mathcal{L} \otimes E)$ for suitable rank two $E$ and any line
  bundle $\mathcal{L}$.  Since $c_1 (E) \equiv c_1 (E \otimes
  \mathcal{L}) \ \mod(2)$ and $\ff_k \cong_{\mathrm{diffeo}}
  \ff_{k+2}$, to find the diffeomorphism type of the base from a
  monodromy equation we need only understand the parity of the first
  Chern class of 
  any suitable bundle $E$ above.
Specifically, in the complex case the map $X \rightarrow \sss$ is the
map defined by the sheaf $\omega_{X / \pp^1}$ and we know
$$c_1 (f_* \omega_{X / \pp^1}) = c_1 (\Lambda^2 
f_* \omega_{X / \pp^1}) = c_1 (\lambda);$$  
and more generally the sphere fibres of the ruled surface
are the projectivisations of spaces of holomorphic sections of the
  canonical bundles down the fibres of $X$.  Thus for any fibration
  $X$ we can take $\det E = \lambda$. From (\ref{fracsign}) and
  (\ref{maintheorem}) we know 
$$4c_1(\lambda) -n-s \ = \ -\frac{3}{5}n - \frac{1}{5}s.$$
It follows that the parity of $10c_1(\lambda)$ and the parity
of $n+2s = 10m$ coincide, and that is precisely what we require.
\end{proof}

\noindent We remark for completeness that these branched coverings can be used
to give a classification of complex genus two fibrations without
reducible fibres.  For in these cases the branch locus is a complex
curve and such curves are determined to smooth isotopy by their
connectivity and homology class.  In
 the connected case, we can ``canonically'' choose such a curve as
 follows; write $s_0$ and $s_{\infty}$ for the homology classes of
 sections of a ruled surface $\ff_k$ of square $\pm k$ respectively
 and $F$ for the homology class of a fibre.  For a curve in a class
 $|r s_{0} + mF|$ choose $m$ fibres and $r$ sections\footnote{For the
   covers to be genus two fibrations we will need $r=6$ since a genus
   two curve double covers a sphere over six points; for the double
   cover to exist the branch divisor, and hence $m$, will have to be
   even.} meeting in disjoint nodes,
 giving a curve of $r+m$ components.  Now under a complex deformation
 of the fibration the branch locus will be perturbed to a neighbouring
 smooth complex curve.  We can clearly arrange for all the fibres to
 lie over some small complex disc in the base $\pp^1$ and for all the
 nodes from intersection points of sections to lie outside the
 preimage of this disc.  The monodromy of the Lefschetz fibration
 arises entirely from a neighbourhood of the $m$ fibres in the branch
 locus and the other nodes, since these are the only places where we
 obtain singular fibres upstairs.  Smoothing the nodal branch locus
 over suitably chosen disjoint discs in the base sequentially, we see
 that the Lefschetz  fibrations arising from
 different branch loci coincide according to
\begin{Eqn}
\begin{array}{rcl}
| r s_{0} + (2m) F|_{\ff_k} & \sim & |r s_{0}|_{\ff_k}
\ \#_{\mathrm{fibre}} \ | rs_{0}  + 2mF |_{\ff_0} \\
& \sim & |rs_{0}|_{\ff_k} \ \#_{m \ \mathrm{fibres}} \ |rs_{0} +
2F|_{\ff_0} \\ & \sim & \#_{i=1}^k \ |rs_{0}|^{(i)}_{\ff_1} 
\ \#_{j=1}^m \ |rs_{0} +2F|^{(j)}_{\ff_0}.
\end{array}
\end{Eqn}
Thus all of the fibrations arising from connected branch loci can be
expressed as fibre sums of two basic pieces.
We have not yet considered disconnected branch loci.  Suppose then
the branch locus is a curve in the class $|as_0| + |bs_{\infty}|$
comprising two disjoint smooth components.  From standard results on
rational surfaces \cite{Hartshorne} it follows that $b=1$.  In
order for the class to contain six sections (and hence yield a genus
two fibration) we are reduced to the case $|5s_0| + |s_{\infty}|$ on
$\ff_{k}$.  Moreover we can write $k=2l$ to obtain an even divisor (as
an element of the Picard group, necessary for the existence of the
double cover).  Again by standard results this class contains no smooth
connected curve, and by the fibre summation trick, it is enough to
understand the fibration when $l=1$;  note that we can build a
disconnected branch locus from only disconnected pieces, and all the
pieces in a fibre sum decomposition will be of this special form
$|5s_0 + s_{\infty}|$ in some $\ff_{2m}$.

\noindent We have reduced all complex genus two fibrations with no
reducible fibres to fibre sums of three basic examples.  These can be
computed in a variety of ways and we obtain the following classification result
(originally due to Chakiris \cite{Chakiris} by similar methods and
discovered independently if much later by the present
author\footnote{Chakiris' result has appeared in numerous articles but
  never accompanied by any kind of proof; his original work, in
  part being more ambitious, is
  somewhat dense.  It therefore seems reasonable to present a version
  of the argument here.}):

\begin{Thm}{\rm(Chakiris)}\label{classify}\qua 
Assume that a genus two fibration has no reducible fibres and K\"ahler
  total space.  Then it is a fibre sum of the shape $A^m B^n = 1$
  or $C^p = 1$ where $m,n,p \in \zz_{\geq 0}$ and the basic words
  $A,B,C$ are given by:
\begin{Eqn}
\begin{array}{rcl}
A\co  \ (\delta_1 \delta_2 \delta_3 \delta_4 \delta_5 \delta_5 \delta_4
\delta_3 \delta_2 \delta_1)^2 = 1 & \sim & X
\stackrel{2:1}{\longrightarrow} \ff_0 \supset |6s_{\infty} +2F| \\
B\co  \ (\delta_1 \delta_2 \delta_3 \delta_4 \delta_5)^6 = 1 & \sim & X
\stackrel{2:1}{\longrightarrow} \ff_1 \supset |6s_0 | \\
C\co  \ (\delta_1 \delta_2 \delta_3 \delta_4 )^{10} = 1 & \sim & X
\stackrel{2:1}{\longrightarrow} \ff_2 \supset |5s_0| \amalg
|s_{\infty}|.
\end{array}
\end{Eqn}
\end{Thm}

\noindent It follows that if a genus two K\"ahler Lefschetz fibration
contains no separating vanishing cycles, then the total space is
simply-connected.
The notation $X \stackrel{d:1}{\longrightarrow} Y \supset B$
indicates 
that $X$ is a $d$--fold branched cover of $Y$ totally ramified along
the locus $B$.  The Dehn twists $\delta_i$ are about curves in the
``standard position'' on a Riemann surface (c.f. \cite{Birman} for
instance).  Note also that for the second family of words $C^p = 1$
the monodromy group of the fibration is not full; but to obtain an
exhaustive list it is sufficient to take all fibre sums by the
identity diffeomorphism.  The word $B=1$ corresponds to the genus two
pencil on a $K3$ which double covers $\pp^2$ branched over a sextic.


\section{Further applications}

At higher genera the moduli spaces are not affine, and
signature is \emph{not} in general determined by the combinatorial
equivalence class of the fibration.  However, for fibrations by
hyperelliptic curves the statement is much easier.  Recall that in the
$(3g-3)$--dimensional moduli space $M_g$ there is a distinguished
$(2g-1)$--dimensional locus of \emph{hyperelliptic curves}, which forms
a complex analytic subvariety $\curly{H}_g$.  We can take the closure
of this subvariety in $\overline{M}_g$.  
The cohomology classes $c_1 (\lambda), \delta_0, \ldots,
\delta_{[g/2]}$ form a basis for $H^2 (\overline{M}_g)$; here
$\delta_0$ is the irreducible component of $\overline{M}_g \backslash
M_g$ corresponding to generically irreducible curves, and the
$\delta_i$ are the components corresponding to curves which
generically have one component of genus $i$ and the other of genus $g-i$.
We can restrict these classes to $\overline{\curly{H}}_g$ where we
denote them by the same symbols.  The following relation seems to be
new:

\begin{Lem}
For any hyperelliptic Lefschetz fibration there is an inequality
$$(8g+4) c_1 (\lambda) \ \geq \ g \delta_0 \ + \ \sum_{h=1}^{[g/2]}
4h(g-h)\delta_h.$$
which is an equality when the base is a two-sphere.
\end{Lem}

\noindent But this follows immediately from (\ref{maintheorem}) and a
fractional signature formula (generalising the one for genus two
above) which has been given by Endo \cite{Endo}.

\noindent As a more significant application we 
provide the answer to a conjecture of Amoros, Bogomolov, Katzarkov and Pantev
\cite{ABKP}.  The monodromy group of a Lefschetz fibration is the subgroup
of the mapping class group generated by the Dehn twists about the
vanishing cycles; that is, the image of the representation
$\pi_1(\pp^1 \backslash \{ \Crit \}) \rightarrow \Gamma_g$.  Recall
that the Torelli group is the subgroup of the mapping class group
comprising elements which act trivially on homology; thus the
monodromy group of $X \rightarrow \pp^1$ is contained in the Torelli
group if and only if all the vanishing cycles are separating.  (We
will sometimes refer to such a fibration as a ``Torelli fibration''.)

For K\"ahler fibrations this can clearly never happen; here is one
classical proof.  By taking the Jacobians down the fibres, a
K\"ahler  fibration gives a map from the smooth locus $\pp^1
\backslash \{ \Crit \}$ to the moduli space of principally polarised
abelian varieties $\mathcal{A}_g$.  This map is holomorphic, and if
the monodromy is trivial on homology groups and hence Jacobians, we
can canonically complete to a holomorphic map of the \emph{closed
  sphere} into $\mathcal{A}_g$.  But this is possible only if the map
is constant, for $\mathcal{A}_g$ is well known to be a bounded complex
domain, and in particular contains no non-trivial holomorphic spheres.

\noindent In \cite{ABKP} it is conjectured\footnote{The conjecture appeared in
  an early draft of the paper; the proof given here was incorporated
  as an Appendix to a later draft at the request of the authors.} that
  in fact there 
  can be no fibrations 
with monodromy group contained in the Torelli group even in the
absence of the K\"ahler assumption (used via holomorphicity of the map
to $\mathcal{A}_g$ above).  Indeed this is true:

\begin{Thm}
There are no Lefschetz fibrations with monodromy group a subgroup of
the Torelli group.
\end{Thm}

\begin{proof}
Suppose such an $(X,f)$ exists; we work by contradiction.
We mentioned before that another approach to signature uses Wall's
non-additivity; an easy case of this approach \cite{Ozbagci} shows
that if a fibration has $\delta$  separating vanishing cycles and no
others, then its signature is given by $-\delta$.  Comparing to the
formula (\ref{maintheorem}) this forces $\langle c_1(\lambda),
[\sss^2]\rangle = 0$.  We will prove that this quantity is positive.

We have $c_2 (X) = 4-4g+\delta = 2-2b_1(X)+b_2(X)$ for any Lefschetz
fibration with $\delta$ singular fibres; in our case, since all the
vanishing cycles are homologically trivial, $b_1(X) = 2g$ and hence
$b_2(X) - 2 = \delta$.  Now the number of disjoint exceptional
$(-1)$--spheres in $X$ is bounded by $b_2(X)$ since each contributes a
new homology class.  The fibre of the fibration has
self-intersection zero; fibre-summing $X$ with itself by the identity
if necessary, we can see that if any fibration with monodromy group
contained in the Torelli group exists, then one exists for which a
section also has self-intersection even and in particular is not
exceptional.  Thus the number of $(-1)$--spheres in $X$ may be assumed
to be bounded above by $b_2(X) - 2 = \delta$.

We invoke a powerful theorem of Liu:  for any minimal symplectic
4--manifold $Z$ which is not irrational ruled, $c_1^2 (Z) \geq 0$.
In particular, blowing down all $(-1)$--spheres in $X$, we see that 
$$c_1^2(X) = c_1^2 (X_{\mathrm{min}}) - \# \{ (-1) \mathrm{-spheres}
\} \geq -\delta$$
and hence 
$$c_1^2(X) + c_2(X) \geq 4-4g$$
provided $X$ is not irrational ruled.  But we also know
$$\frac{c_1^2(X) + c_2 (X)}{12} = \frac{\sigma(X) + \chi(X)}{4} =
1-g+\langle c_1(\lambda), [\sss^2_X] \rangle.$$
\noindent Rearranging this gives $\langle c_1(\lambda), [\sss^2]
\rangle > 0$ which is the required contradiction.  That leaves only
the case of $X$ irrational ruled, but simple topological conditions
show that such a manifold can have no Torelli fibration.  For any
Torelli Lefschetz fibration of $X = \Sigma_h \tilde{\times} \sss^2$ has
fibres of genus $h$ by considering $H_1(X)$, and then looking at
$\pi_1(X)$ we see that in fact all the vanishing cycles must be
nullhomotopic and hence $X$ was a trivial product. 
\end{proof}

\noindent There are fibrations over the torus with no
critical fibres which have monodromy group contained in the Torelli
group; just take $\sss^1 \times Y$ for $Y$ the mapping torus of a Dehn
twist about a separating curve.  The Lefschetz fibrations of these
manifolds, however, cannot be Torelli.  

\begin{Cor}
For any Lefschetz fibration $X$, $\sigma(X) + \delta > 0$.  Moreover the
sphere $\sss^2 _X \subset \mgbar$ defined to isotopy by $X$ has
positive ``symplectic volume''; the symplectic Weil--Petersson form on
$\mgbar$ evaluates positively on $\sss^2 _X$.
\end{Cor}

\noindent The corollary again follows by comparing to Ozbagci's work; for he
proves that each critical fibre changes the signature by $0, \pm 1$
as you add handles to a trivial bundle over the disc.  Since we now
know there must be a non-separating vanishing cycle, and his
construction allows us to add that handle first and change the
signature by a non-negative amount, we obtain the first result.  In particular,
it follows that for \emph{any} smooth sphere in $\overline{M}_g$ with
correct geometric intersection behaviour at the compactification
divisors, the Chern class $c_1 (\lambda)$ evaluates
positively.  The K\"ahler Weil--Petersson form is not quite a
positive rational multiple of this Chern class and the boundary
divisors; rather we have the identity (technically on the moduli stack)
$$\frac{1}{2\pi^2} \omega_{WP} \ = \ 12c_1 (\lambda) - \delta.$$
An easy computation, however, shows that $c_1 (\lambda) > 0$
and $c_1 ^2 (X) + 8(g-1) > 0$ together imply that the Weil--Petersson
  form evaluates positively.  This is made more interesting by the
  following remark. Wolpert has computed the intersection ring of
$\overline{M}_g$ and it follows from his computation that this is
\emph{not} a purely homological statement \cite{Wolpert}; there are
  two dimensional
homology classes which have positive intersections with all the
components of $\overline{M}_g \backslash M_g$ but not with $c_1
(\lambda)$.  Thus the corollary is a kind of ``symplectic ampleness''
phenomenon, reliant on the \emph{local} positivity of intersections with the
stable locus.
\vfil \eject



\begin{thebibliography}

\bibitem{ABKP}
{\bf J Amoros et\,al},
 {\it Symplectic {L}efschetz fibrations with arbitrary
  fundamental groups},
 preprint (1998) with an Appendix by Ivan Smith

\bibitem{arabel}
{\bf E Arabello}, {\bf M Cornalba},
 {\it The {P}icard groups of the moduli spaces of
  curves},
 Topology {26} (1987) 153--71

\bibitem{Atiyaheta}
{\bf M Atiyah},
 {\it The logarithm of the {D}edekind $\eta$--function}, Math. Ann.
  {278} (1987) 335--380 

\bibitem{APS}
{\bf M Atiyah et\,al},
 {\it Spectral asymmetry and {R}iemannian geometry {I}},
  Math. Proc. Camb. Phil. Soc. {77} (1975)  43--69 

\bibitem{BPV}
{\bf W\,Barth}, {\bf C\,Peters}, {\bf A\,Van\,de\,Ven},
 {\it Compact complex surfaces},
  Springer (1984)

\bibitem{Birman}
{\bf J Birman},
 {\it Braids, links and mapping class groups}, 
Princeton University
  Press (1975)

\bibitem{Brown}
{\bf K\,S Brown},
 {\it Cohomology of groups}, Springer (1982) 

\bibitem{Chakiris}
{\bf K\,N Chakiris},
 {\it The monodromy of genus two pencils},
 Ph\,D  thesis
  Colombia University (1978)

\bibitem{Endo}
{\bf H Endo},
 {\it Meyer's signature cocycle and hyperelliptic fibrations},
  preprint (1998) 

\bibitem{Freed}
{\bf D Freed},
 {\it On determinant line bundles}, 
from: ``Mathematical Aspects of String
  Theory {I}'', S\,T Yau (editor),  World Scientific (1987)

\bibitem{Fuller}
{\bf T Fuller}, 
{\it Genus two {Lefschetz} fibrations}, preprint (1997) 

\bibitem{Hartshorne}
{\bf R Hartshorne},
 {\it Algebraic geometry},
 Springer (1977)

\bibitem{Matsumoto}
{\bf Y Matsumoto},
 {\it Lefschetz fibrations of genus two - a topological
  approach},
 from: ``The 37th Taniguchi Symposium on topology and Teichm\"uller spaces'',
  S Kojima et\,al  (editors),  World Scientific (1996) 

\bibitem{MumfordBirk}
{\bf D Mumford},
 {\it Towards an enumerative geometry of the moduli space of
  curves},
 from: ``Arithmetic and Geometry'', A Artin and J Tate (editors),  vol~II
  Birkh\"auser (1983)

\bibitem{Ozbagci}
{\bf B Ozbagci},
 {\it Signatures of {L}efschetz fibrations},
 preprint (1998) 

\bibitem{ST}
{\bf B\,Siebert}, {\bf G\,Tian},
 {\it On hyperelliptic $C^{\infty}$--Lefschetz fibrations
  of 4--manifolds},
 preprint (1999)

\bibitem{ivanthesis}
{\bf I Smith},
 {\it Symplectic geometry of Lefschetz fibrations},
 Ph\,D  thesis,  Oxford University (1998)

\bibitem{Wolpert}
{\bf S Wolpert},
 {\it On the homology of the moduli space of stable curves}, 
Ann.  of Math. {118} (1983)  491--523 

\bibitem{Zariski}
{\bf S Zariski},
 {\it Algebraic surfaces},
 2nd ed. Springer (1977) with appendices
  by {D~Mumford} 

\end{thebibliography}
\end{document}